# An efficient online successive reanalysis method for dynamic topology optimization

Shuhao Li[a,b] Hu Wang[a,b*], Jichao Yin[a,b], Daozhen Guo[a,b], Guangyao Li[b]

[a]*State Key Laboratory of Advanced Design and Manufacturing for Vehicle Body, Hunan University, Changsha 410082, People's Republic of China;*

[b]*Beijing Institute of Technology Shenzhen Automotive Research Institute, Shenzhen 518000, People's Republic of China*

[*]Hu Wang: wanghu@szari.ac.cn

**Highlights**

- An efficient online successive dynamic reanalysis method is proposed
- A POD-based approximate dynamic displacement strategy is proposed
- Update and convergence criteria of equivalent topology are given
- Several large-scale cases demonstrate the effectiveness of the proposed method

**Abstract**

In this study, an efficient reanalysis strategy for dynamic topology optimization is proposed. Compared with other related studies, an online successive dynamic reanalysis method and POD-based approximate dynamic displacement strategy are integrated. In dynamic reanalysis, the storage of the stiffness matrix decomposition can be avoided and the reduced basis vectors should be updated successively according to the structural status in each iteration. Therefore, the bottleneck of combined approximation method for large-scale dynamic topology optimization can be handled. Sequentially, the Proper Orthogonal Decomposition (POD) is employed to obtain the approximate dynamic displacement, in which the Proper Orthogonal Mode (POM) of the displacement field is employed to establish the approximated equivalent static loads of Equivalent Static Load (ESL) method. Compared with the exact equivalent static loads at all the time intervals, the number of equivalent static loads is significantly reduced. Finally, the 2D and 3D test results indicate that the proposed method has remarkable speed-up effect on the premise of small relative error, support the strength of the proposed strategy.

**Keywords** Equivalent static loads; Proper orthogonal decomposition; Combined

approximation method; Dynamic topology optimization

## 1. Introduction

In recent years, topology optimization methods have developed rapidly and have been widely used in engineering applications. Using topology optimization, engineers can quickly obtain the optimal structural design of products in the early stage of design [1]. Most of the engineering cases are dynamic problems, which need to consider the time-dependent characteristics. Generally, the research of dynamic topology optimization method can be classified as two branches: one is the topology optimization of frequency-domain dynamic problems, and the other is time-domain dynamic problems. For the topology optimization of frequency-domain problems, Díaaz and Kikuchi [2] used the homogenization method to maximize the natural frequency of the structure, which might be the earliest research on dynamic topology optimization. Subsequently, a large number of topology optimization studies based on eigenvalue problem emerged [3-5]. For the topology optimization of time-domain problems, the structure is excited by external load, and the displacement, velocity and acceleration of the structure are time-dependent quantities. Considering these dynamic responses as the objective functions, the topology optimization of the structure is studied [6, 7]. In this study, the topology optimization of time-domain problems is discussed.

For most of related works, the complete analysis is employed. In term of accuracy, the complete analysis-based solution can be regarded as high-fidelity simulation. However, simulations of these high-fidelity models, especially for time-domain problems, entail expensive computations. The repeated calculation in each topology optimization cycle further increases the computational cost, rendering direct numerical simulation infeasible. In order to overcome this bottleneck, the model order reduction has been applied to the dynamic analysis process [8, 9]. Model order reduction technology reduces the computational cost by projecting the original problem into a low dimensional reduced space, and realizes the order reduction of large-scale problems [10]. It has been widely used in real-time simulation[11, 12], model-based optimization [8] and other research fields[13, 14].

In these model order reduction methods, a kind is to construct the approximate solution of the problem through a series expansion. The parameters of the expansion term are determined by the solution of the system constructed by several basis vectors. Among these model order reduction techniques, reanalysis method has been widely used in the field of finite element simulation. Reanalysis is a fast calculation method to quickly obtain the response of the modified structure through the information of the initial structure, including the response of the initial structure and the intermediate variables. The reanalysis method avoids the repeat calculation of the response of the modified structure. By construction of reduced space based on the initial structure information, the response of the modified structure is obtained by approximate calculation, which significantly reduces the computational burden. Among many reanalysis methods, Combined Approximation (CA) method has both the efficiency of local approximation and the accuracy of global approximation, so it might be the most popular one. The CA method was first proposed by Kirsch [15]. After continuous improvement in recent decades, it has been successfully applied to the fast calculation of static problems and frequency-domain dynamic problems [16-18]. However, the CA method is rarely used in time-domain dynamic problems, while it is also difficult for the CA to handle large modification cases. To solve the large-scale time-domain dynamical problems, an efficient reanalysis strategy in the FE framework is proposed, which is named Online Successive Dynamic Combination Approximation (OSDCA) method. The core idea of this method is that the updated structure of each dynamic optimization iteration is based on the structure information of previous adjacent iteration. The topology changes little between adjacent optimization iterations, which makes the application of reanalysis method obtain better approximation effect. In each optimization iteration, the reduced basis vector is constructed online to approximate the updated structure, which avoids the problem of excessive memory consumption of matrix decomposition when solving large-scale equilibrium equations. Thus, the dynamic response can be approximated quickly and accurately in each topology iteration.

Based on the approximate dynamic response obtained by the above OSDCA method, the corresponding topology optimization can be constructed. Among dynamic topology optimization methods, Equivalent Static Load (ESL) method, proposed by Choi and Park [19], as an efficient simplified method to deal with dynamic problems, has been widely used. Its core idea is to transform the dynamic topology optimization problem into a static topology optimization problem under multiple external loads. Based on the dynamic displacement, a series of equivalent static loads, namely ESL set, are constructed. The static displacement generated by these loads is consistent with the displacement corresponding to arbitrary time interval. Then, the static topology optimization process under equivalent static loads is executed. The calculation time of static topology optimization process is directly proportional to the number of ESLs. Since the dynamic analysis process should be divided into a large number of time intervals, it is not advisable to construct the exact ESLs. It is necessary to select the displacement solution covering the main information to construct the approximated ESLs. The ESL set can be calculated at some peaks of important locations [20, 21] .However, this method is limited to the location of pre-assumption and should have some locality[22]. It is still difficult for the present ESL to apply for time-dependent problem because the characteristics of the entire time domain are difficult to be extracted.

To overcome this issue, POD, as a data compression technology, which can extract important information in the data set, is introduced. Actually, the POD has been successfully applied to many aspects of scientific research and engineering, including model analysis [23], shape optimization [24] and topology optimization [25, 26]. In these applications, POD is mostly used for the reduced order model construction process in the dynamic analysis stage. Lee and Cho [27] propose to apply POD in the static optimization stage of the ESL-based optimization framework to directly obtain approximated ESLs, which has been successfully applied to size optimization. For these reasons, the POD method might be the suitable tool for the ESL-based topology optimization. In this proposed POD-based approximated ESL method, POD is used to comprise the essential information contained in dynamic displacement field. Then the approximated ESLs are calculated based on displacement's POD approximation. Since the snapshot comes from the displacement solution of all time intervals, the exact ESL set is reduced globally. Thus, several ESLs can represent the exact ESL set at all the time intervals, which greatly improves the calculation efficiency.

The rest of this paper' organization is as follows. Section 2 introduces the dynamic topology optimization problem. Section 3 develops the OSDCA method, which is applied to the dynamic simulation based on Newmark method. Section 4 introduces the POD for dynamic displacement field, and the approximated ESLs are calculated. Based on the static displacement under these approximated ESLs, the objective function of dynamic topology optimization problem is constructed. Section 5 provides the results of the three numerical examples, including 2D and 3D beam structures. Finally, the conclusions and suggestions are given.

## 2. Dynamic structural topology optimization problem

The purpose of topology optimization is to optimally distribute material within a design domain. The focus of this study is to improve the efficiency of topology optimization problem of large-scale dynamic system. With the increase of the scale of the problem, the time consumption of the dynamic analysis increases significantly. When applied to the topology optimization, repeated analysis in each iteration lead to low efficiency of the whole optimization process, and even render optimization infeasible. In this section, the topology optimization problem with linear elastic structures subjected to dynamic loading is given. A reduction strategy is put forward to improve the efficiency of dynamic system topology optimization.

### 2.1. Dynamic equilibrium equations

The governing dynamic equilibrium equations for a multi-DOF can be expressed as follows:

$$\boldsymbol{M}\ddot{\boldsymbol{d}}(t) + \boldsymbol{C}\dot{\boldsymbol{d}}(t) + \boldsymbol{K}\boldsymbol{d}(t) = \boldsymbol{F}(t) \tag{1}$$

where $\boldsymbol{M}$, $\boldsymbol{C}$, $\boldsymbol{K}$ denote mass, damping, stiffness matrices respectively. Displacement, velocity, acceleration and external force vectors are functions of the time variable $t$.

### 2.2. Dynamic topology optimization strategy

In ESL-based dynamic topology optimization, the dynamic problem should be transformed into static problem under a multiple loading condition. At each dynamic optimization iteration, dynamic analysis and static analysis under multiple ESLs account for the main time. For large-scale problems, precise dynamic analysis and employ of exact ESLs should be achieved improbably. In this study, efficient approximation methods are developed for these two aspects. The framework of the suggested method is shown in Figure. 1.

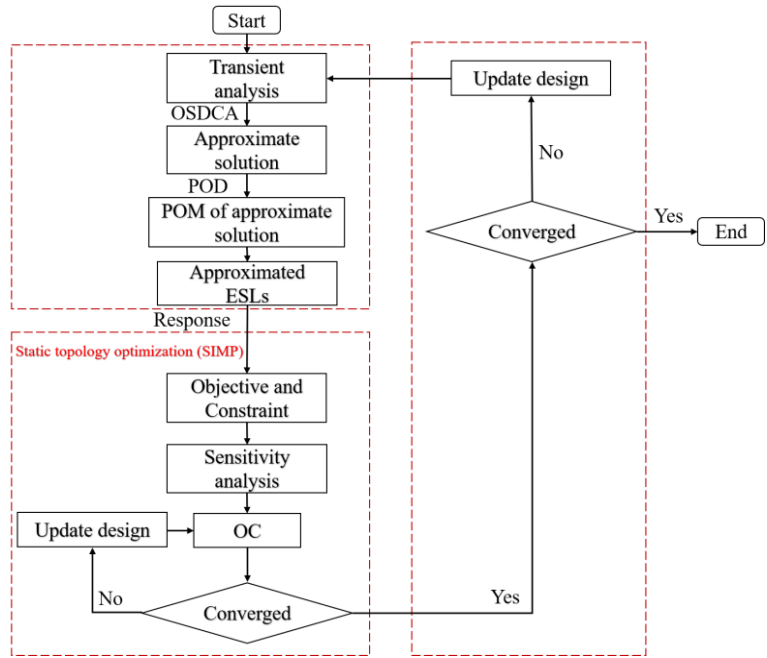

Fig. 1 Framework of efficient dynamic topology optimization

The suggested framework is composed of two stages. In the first stage, dynamic analysis response is obtained by using the reduction strategy proposed in Section 3. Firstly, the dynamic equilibrium Eq. (1) is solved. In OSDCA method, the original dynamic problem should be represented by several basis vectors in the reduced space. Then, by introducing POD method, the main characteristics of dynamic displacement field are extracted, that is, the POM of dynamic displacement snapshot. Finally, according to the obtained displacement's POD approximation, the approximated ESLs should be achieved. Thus, the dynamic problem can be transformed into a static problem under the several actions of approximated ESLs, and the static displacement response is obtained finally.

The second stage is to solve the topology optimization problem. Solid Isotropic Material with Penalization (SIMP) method based on Heaviside projection function is employed. After finding the optimal solution of the static topology optimization, a transient analysis is executed by using the optimum to update the ESLs. Then, the next static optimization starts with an updated multiple loading condition (updated ESLs). As a consequence, the converged solution of the outer loop becomes the final optimum of the topology optimization for dynamics. The specific details will be introduced in details in the following sections.

## 3. Dynamic structural reanalysis method

As the scale of the problem increases, the time-consuming proportion of the dynamic analysis in the entire dynamic optimization should be correspondingly increased. In this section, the OSDCA method is proposed to improve the efficiency of dynamic analysis.

### 3.1. Newmark integration method for dynamic system

The dynamic process is solved by Newmark method, which is an implicit numerical integration method for solving partial differential equations. Because of its stability, it is widely used to solve problems in dynamic field. Different numerical values of parameters in Newmark method ($\alpha$ and $\beta$) determine different numerical integration schemes. Here, $\alpha$ and $\beta$ take 0.25 and 0.5 respectively. At this time, Newmark method is equivalent to constant average acceleration method, which is an unconditionally stable integration scheme, that is, the size of time interval $\Delta t$ does not affect the stability of the solution, which makes the method suitable for the dynamic problems.

According to Eq. (1), solved by the Newmark method, the recursive formula with $\boldsymbol{d}(t + \Delta t)$ as the unknown quantity is adopted. The main solution steps are as follows:

Firstly, the finite difference relation is substituted into Eq. (1) to obtain the coefficient of the integration:

$$c_0 = \frac{1}{\alpha \Delta t^2}, c_1 = \frac{\beta}{\alpha \Delta t}, c_2 = \frac{1}{\alpha \Delta t}, c_3 = \frac{1}{2\alpha} - 1$$
$$c_4 = \frac{\beta}{\alpha} - 1, c_5 = \Delta t \left(\frac{\beta}{2\alpha} - 1\right), c_6 = \Delta t(1 - \beta), c_7 = \beta \Delta t \tag{2}$$

Calculating the effective stiffness matrix $\widehat{\boldsymbol{K}}$:

$$\widehat{\boldsymbol{K}} = \boldsymbol{K} + c_0 \boldsymbol{M} + c_1 \boldsymbol{C} \tag{3}$$

Obtaining the decomposed form:

$$\widehat{\boldsymbol{K}} = \boldsymbol{LDL}^T \tag{4}$$

Then, for each time interval ($t = 0, \Delta t, 2\Delta t ...$), calculating the effective load vector $\widehat{\boldsymbol{F}}(t + \Delta t)$:

$$\begin{aligned}\widehat{\boldsymbol{F}}(t + \Delta t) = \boldsymbol{F}(t + \Delta t) &+ \boldsymbol{M}\left(c_0 \boldsymbol{d}(t) + c_2 \dot{\boldsymbol{d}}(t) + c_3 \ddot{\boldsymbol{d}}(t)\right) \\ &+ \boldsymbol{C}(c_1 \boldsymbol{d}(t) + c_4 \dot{\boldsymbol{d}}(t) + c_5 \ddot{\boldsymbol{d}}(t))\end{aligned} \tag{5}$$

Calculating the nodal displacement $\boldsymbol{d}(t + \Delta t)$ in each time interval:

$$\boldsymbol{LDL}^T \boldsymbol{d}(t + \Delta t) = \widehat{\boldsymbol{F}}(t + \Delta t) \tag{6}$$

The nodal acceleration $\ddot{\boldsymbol{d}}(t + \Delta t)$, nodal velocity $\dot{\boldsymbol{d}}(t + \Delta t)$:

$$\ddot{\boldsymbol{d}}(t + \Delta t) = c_0 \big(\boldsymbol{d}(t + \Delta t) - \boldsymbol{d}(t)\big) - c_2 \dot{\boldsymbol{d}}(t) - c_3 \ddot{\boldsymbol{d}}(t) \tag{7}$$

$$\dot{\boldsymbol{d}}(t + \Delta t) = \dot{\boldsymbol{d}}(t) + c_6 \ddot{\boldsymbol{d}}(t) + c_7 \ddot{\boldsymbol{d}}(t + \Delta t) \tag{8}$$

From the above calculation steps, it can be seen that the most time-consuming part of the dynamic analysis is the solution of equilibrium equations in Eq. (6). The Section 3.2 will introduce how to replace these equilibrium equations in the original problem with the reduced equilibrium equations.

### *3.2. Dynamical reanalysis formulation*

According to above mentioned issue, a dynamical reanalysis is proposed in this section to improve the efficiency of Newmark method. The initial problem considered in this study can be formulated as (damping matrix $\boldsymbol{C}$ is omitted to facilitate the derivations):

$$\boldsymbol{M}_0 \ddot{\boldsymbol{d}}(t) + \boldsymbol{K}_0 \boldsymbol{d}(t) = \boldsymbol{F}(t) \tag{9}$$

In dynamic topology optimization, the material distribution of the structure should be updated in each iteration, and the corresponding structural changes are reflected in the equation as $\Delta \boldsymbol{K}$ and $\Delta \boldsymbol{M}$, so that the modified $\boldsymbol{K}$ and $\boldsymbol{M}$ are given by

$$\boldsymbol{K} = \boldsymbol{K}_0 + \Delta \boldsymbol{K}$$
$$\boldsymbol{M} = \boldsymbol{M}_0 + \Delta \boldsymbol{M} \tag{10}$$

The equilibrium equation of the updated structure can be expressed as Eq. (11). The purpose of dynamic reanalysis is to avoid solving the equilibrium equation of the updated structure and obtain the approximate response of the modified structure.

$$\boldsymbol{M}\ddot{\boldsymbol{d}}(t) + \boldsymbol{K}\boldsymbol{d}(t) = \boldsymbol{F}(t) \tag{11}$$

Substituting the modified $\boldsymbol{K}$ and $\boldsymbol{M}$ into Eq. (3) and Eq. (5), equation Eq. (6) can be rewritten as a static reanalysis form. Thus, a dynamical reanalysis problem is rewritten as the reanalysis form in each time interval.

$$(\widehat{\boldsymbol{K}}_0 + \Delta \widehat{\boldsymbol{K}})\boldsymbol{d}(t + \Delta t) = \widehat{\boldsymbol{F}}(t + \Delta t) \tag{12}$$

where

$$\Delta \widehat{\boldsymbol{K}} = \Delta \boldsymbol{K} + c_0 \Delta \boldsymbol{M} \tag{13}$$

### *3.3. Online successive dynamic combined approximation method*

In the process of dynamic topology optimization, the equilibrium equation needs to be solved in each iteration even if the change of structure is small, so the amount of calculation will increase exponentially.

As mentioned in the introduction, reanalysis method can quickly obtain the response after structural change without the repeat analysis, and CA method is widely used because of its global approximation and efficiency. The details of the classical CA method can be found in Akgün et al [28]. Based on the classical CA method, a dynamic topology optimization method is proposed. Actually, time-dependent dynamic reanalysis methods suggested by Gao [29] proposed a feasible Adaptive Time-based Global Reanalysis (ATGR) algorithm for Newmark-$\beta$ method. However, the scale of their numerical examples is limited and the adaptive strategy involves too many operations. Therefore, it is difficult for the ATGR to solve large scale problems and the efficiency improvement of the ATGR is not good as expected.

There are two difficulties in applying reanalysis to dynamic topology optimization. Firstly, the structure changes greatly after multiple dynamic topology iterations or even two iterations, and it can be regarded as a global change. In order to obtain an accurate approximation result, the reduced physical matrix, such as stiffness matrix should be updated at each iteration. Secondly, the construction of basis vectors involves multiplier operations, such as solution of equilibrium equations, matrix decomposition and others. This number of operation numbers and the storage of the initial calculation information should limit the feasibility of dynamic reanalysis. To overcome the above bottlenecks, the online successive dynamic combination approximation method is proposed and the flowchart of this method is presented in Figure 2.

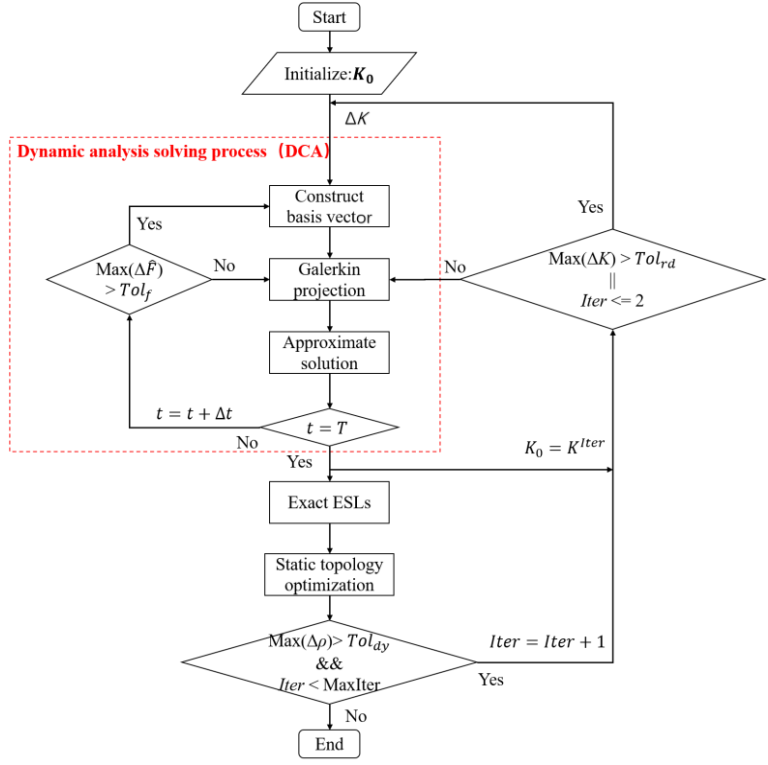

Fig. 2 The flowchart of the OSDCA

Firstly, the stiffness matrix $K_0$ is initialized, and the variation of structural stiffness matrix $\Delta K$ in the dynamic topology optimization iteration is calculated, so as to construct the basis vectors. The displacement response of each time interval is obtained by Galerkin projection. Besides, the basis vectors are also adaptively reconstructed in each dynamic analysis according to the given threshold $Tol_f$. Then, the design variable $\rho$ is updated through the static topology optimization process. In each iteration of dynamic topology optimization, the structure $K_0$ is updated to judge whether $\Delta K$ exceeds the tolerance $Tol_{rd}$ (defined in Section 5.1). If it exceeds the tolerance, the basis vectors are reconstructed through the $\Delta K$ between the current iteration structure and the previous optimized structure. If it meets the tolerance requirements, the approximate displacement is obtained directly through Galerkin projection. The core of this method is to constantly update the initial structure information $K_0$ to construct the basis vectors. The details of the approximate displacement calculation by the proposed dynamic reanalysis method are shown in Figure 3.

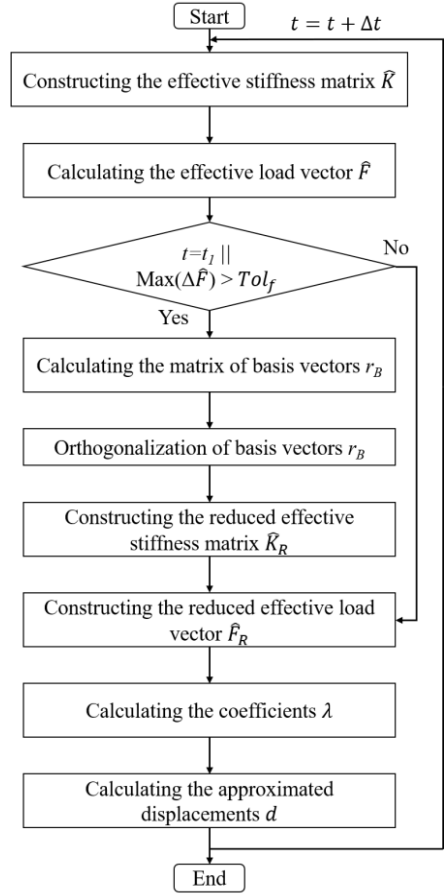

Fig. 3 The flowchart of the DCA

Calculation of $d(t)$ by online dynamic combination approximation method involves the following steps:

(1) Suppose the structural response of each time interval is approximately expressed as

$$\boldsymbol{d} = \lambda_1 \boldsymbol{r}_1 + \lambda_2 \boldsymbol{r}_2 + \cdots + \lambda_s \boldsymbol{r}_s = \boldsymbol{r}_B \boldsymbol{\lambda} \qquad (14)$$

where $s$ is the number of basis vectors, and $s$ is much smaller than the number of DOFs $h$, $\lambda$ the coefficients matrix

$$\lambda = [\lambda_1 \ \lambda_2 \ \dots \ \lambda_s]^T \qquad (15)$$

$\boldsymbol{r}_B$ the basis vectors

$$\boldsymbol{r}_B = [\boldsymbol{r}_1 \ \boldsymbol{r}_2 \ \dots \ \boldsymbol{r}_s] \qquad (16)$$

(2) Construct the initial effective stiffness matrix $\widehat{\boldsymbol{K}}_0$ by Eq. (3) and the modified effective stiffness matrix $\widehat{\boldsymbol{K}}$ is

$$\widehat{K} = \widehat{K}_0 + \Delta\widehat{K} \tag{17}$$

(3) In the 1[th] time interval $t_1$, the reduced basis vectors are constructed as follows:

a) Calculate the effective load vector in the $t_1$ time interval

$$\widehat{F}(t_i) = F(t_i) + M\left(c_0 d(t_{i-1}) + c_2 \dot{d}(t_{i-1}) + c_3 \ddot{d}(t_{i-1})\right), i = 1 \tag{18}$$

b) Calculate the matrix of basis vectors $r_B$, and the basis vectors $r_i$ ($i$=1,2,…,$s$) are given by

$$\Delta\widehat{K} = \widehat{K} - \widehat{K}_0 \tag{19}$$

$$r_1 = \widehat{F}(t_1)/\widehat{K}_0 \tag{20}$$

$$r_i = -Br_{i-1}(i = 1,2,\dots,s)$$
$$B = \Delta\widehat{K}/\widehat{K}_0 \tag{21}$$

c) In order to ensure the linear independence of the basis vectors, the generated basis vectors orthogonalized:

$$[r_1\ r_2\ \dots\ r_s] = r_B \Sigma A^T \tag{22}$$

d) Construct the reduced effective stiffness matrix $\widehat{K}_R$ and effective load vector $\widehat{F}_R$ by

$$\widehat{K}_R = r_B^T \widehat{K} r_B \tag{23}$$

$$\widehat{F}_R = r_B^T \widehat{F}(t_i), i = 1 \tag{24}$$

e) Calculate the coefficients $\lambda$ by solving the reduced system

$$\lambda = \widehat{F}_R / \widehat{K}_R \tag{25}$$

f) Calculate the approximated displacements $d(t_1)$ in the $t_1$ time interval

$$d(t_i) = r_B \lambda, i = 1 \tag{26}$$

(4) In the following time intervals, effective load alteration ($\Delta\widehat{F} = \left|\frac{\widehat{F}_i - \widehat{F}_{i-1}}{\widehat{F}_{i-1}}\right|$) is defined.

If $\Delta\widehat{F} > Tol_f$ (a given threshold), reconstructing the basis vectors, otherwise, calculating the approximated displacements $d(t)$ ($t = t_2\ t_3\ \dots\ t_l$) by Eqs. (18), (24-26), where $l$ is the maximum time interval. In this way, the approximate displacement response $\boldsymbol{d}(t)$ ($t = t_1\ t_2\ t_3\ ...\ t_l$) at all the time intervals is obtained.

By comparing Eq. (12) and Eq. (15), the modified effective stiffness $\widehat{\boldsymbol{K}}$ size is $h \times h$, the reduced effective stiffness $\widehat{\boldsymbol{K}}_R$ is $s \times s$, and $s \ll h$, the scale of the original problem is significantly reduced and the calculation efficiency is improved.

Compared with the ATGR, the basis vectors are obtained by the structural information of adjacent iterations dynamically, and auxiliary operations for generation of multiple sample points should be avoided. Although the ATGR uses global approximation, the accuracy of prediction can't be guaranteed yet due to the large modifications, especially for topology optimization. Moreover, consideration of operation number of ATGR, the computational cost can also be improved. In the suggested way, the basis vectors are constantly updated throughout the dynamic topology optimization cycle, that is, solve several equilibrium equations directly. Therefore, the accuracy should be guaranteed compared with the ATGR. Generally, the calculation time of single construction basis vectors is much less than that of the whole dynamic analysis process, so as to realize efficient online and fast calculation. Moreover, because matrix decomposition and initial structure information storage operations are avoided, it can be applied to large-scale dynamic topology problems.

## 4. POD-ESL dynamic topology optimization

In Section 3.3, the approximate displacement obtained by the OSDCA should be employed to generate snapshot matrix. In this section, POD is introduced to extract the main characteristics in the entire time domain and obtain the displacement's POD approximation. The approximated ESLs are constructed by the POM of the dynamic displacement. Thus, the POD-based approximated ESLs strategy is suggested to reduce the number of ESLs, and the original objective function should be updated for topology optimization based on the response solved previously, and the derivation of analytical sensitivities is discussed in the subsequent subsection.

*4.1. POD-based approximated equivalent static load strategy*

Proper orthogonal decomposition is a data compression technology, which can extract important information in the data set. These data sets can come from the numerical simulation results of complex physical fields (such as flow field and thermal field), experimental data or various image data sources. A series of linear bases, called POD bases, can be obtained by using POD. Based on these basis vectors, the reduced solution with minimum projection error , equivalent to maximizing energy in the projection, can be obtained by interpolation in the original system [30]. A problem can often be expressed or approximated by a series of high-fidelity solutions. These high-fidelity solutions are called the solution flow of the problem, from which a series of basis vectors can be extracted to characterize the problem. According to the decline rate of singular value, the number of basis vectors used to approximately represent the original data can be obviously judged. In this study, the snapshot matrix composed of displacement vectors corresponding to a series of time intervals obtained by using OSDCA method to solve the high-fidelity system. The dynamic displacement can be approximately characterized by a small number of basis vectors obtained by POD.

For convenience, the displacement snapshot matrix is directly used to introduce POD for subsequent research. According to the definition in Section 3.3. The core of POD is to approximate the dynamic displacement snapshots $\{\boldsymbol{d}_i\}_{i=1}^{l}$ simultaneously by a single, normalized vector $\boldsymbol{\gamma}$ as well as possible, that is, to solve the optimization problem:

$$\max_{\boldsymbol{\gamma}\in\mathbb{R}^n}\frac{1}{l}\sum_{i=1}^{l}|<\boldsymbol{d}_i,\boldsymbol{\gamma}>|^2 \qquad (27)$$

Subject to $||\boldsymbol{\gamma}||^2 = 1, \boldsymbol{\gamma} \in \mathbb{R}^n$

where $||\boldsymbol{\gamma}|| = \sqrt{<\boldsymbol{\gamma},\boldsymbol{\gamma}>}$.

In Eq. (27), $l$ is the number of snapshots, that is equal to the maximum time interval.

Apply Lagrange functional to Eq. (27), the eigenvalue problem related to the optimization problem can be written as:

$$\frac{1}{m}UU^T\boldsymbol{\gamma} = \varphi\boldsymbol{\gamma} \qquad (28)$$

where $U$ is the displacement matrix, and with the condition

$$||\boldsymbol{\gamma}|| = 1 \qquad (29)$$

The eigenvalue problem in Eq. (28) is solved to obtain the orthogonal vectors $\{\boldsymbol{\varphi}_j\}_{j=1}^{m}$ of the displacement snapshots $\{\boldsymbol{d}_i\}_{i=1}^{l}$, columns of $U$. These orthogonal vectors are called POD basis of rank $m \leq l$. Select the first $m$ vectors $\{\boldsymbol{\varphi}_j\}_{j=1}^{m}$ to approximately characterize the snapshots $\{\boldsymbol{d}_i\}_{i=1}^{l}$, and these basis vectors have the least approximation error among all rank $m$ approximations to the columns of $U$. This method is the so-called method of snapshots.

In POD method, the selection of the number of basis vectors is important. As a snapshot, the displacement obtained by the high-fidelity model can well represent the problem studied. It has a priori criterion to refer to and can ensure the accuracy of approximation as much as possible, but there is no priori rule for the selection of the number of basis vectors. Rather, heuristic criteria can be used as a reference[31]. Selection criteria is based on the ratio of the energy contained in the first $m$ vectors $\{\boldsymbol{\varphi}_j\}_{j=1}^{m}$ selected in the system to total energy, and the ratio is defined as:

$$\varepsilon(m) = \frac{\sum_{j=1}^{m}\boldsymbol{\varphi}_j}{\sum_{j=1}^{l}\boldsymbol{\varphi}_j} \qquad (30)$$

The key of applying ESL method to solve dynamic problems (dynamic topology optimization problems in this study) is to construct ESL set that produce the same displacement solutions as the original dynamic problems. This makes each time interval correspond to an ESL, and the displacement obtained by each ESL is consistent with the displacement field of the corresponding time interval. This principle is compatible for the goal of minimizing compliance in topology optimization. The exact ESL set is constructed as follows:

$$\boldsymbol{K} \times \{\boldsymbol{d}_1\ \boldsymbol{d}_2\ ...\ \boldsymbol{d}_l\} = \boldsymbol{F} = \{\boldsymbol{f}_1\ \boldsymbol{f}_2\ ...\ \boldsymbol{f}_l\} \tag{31}$$

where $\boldsymbol{K}$ is the stiffness matrix; $\boldsymbol{d}$ the displacement field vector corresponding to each time interval; $\boldsymbol{f}$ the vector of ESLs; $\boldsymbol{F}$ the exact ESLs and $l$ the number of loads, which consistent with time intervals. For arbitrary ESLs, the following conditions are met:

$$\begin{aligned} \boldsymbol{f}_1 &= \boldsymbol{K}\boldsymbol{d}_1 \\ \boldsymbol{f}_2 &= \boldsymbol{K}\boldsymbol{d}_2 \\ &\cdots\cdots \\ \boldsymbol{f}_l &= \boldsymbol{K}\boldsymbol{d}_l \end{aligned} \tag{32}$$

Considering the computational efficiency and the complexity of the dynamic topology optimization problem, it is not advisable to select the exact ESLs at all the time intervals. Because in order to pursue the accuracy of displacement response, a large number of time intervals need to be adopted, which leads to a great increase in the computational burden of static topology optimization. To do so, POD is introduced to extract the main characteristics of displacement field. In fact, various researches have been carried out regarding the POD in dynamic problems [27, 32, 33]. In this study, select the approximate displacement obtained by OSDCA method to construct the snapshots, and calculate the POM of the displacement field $\boldsymbol{U}$ in a global manner:

$$\boldsymbol{U} = \{\boldsymbol{d}_1\ \boldsymbol{d}_2\ ...\ \boldsymbol{d}_l\} = \boldsymbol{\Phi}\boldsymbol{\Sigma}\boldsymbol{A}^T \tag{33}$$

where $\boldsymbol{\Phi}$ is the proper orthogonal mode (POM):

$$\boldsymbol{\Phi} = \{\boldsymbol{\varphi}_1 \; \boldsymbol{\varphi}_2 \ldots \boldsymbol{\varphi}_m\}, m \leq l \tag{34}$$

$\Sigma$ is a diagonal matrix of singular values, and $\boldsymbol{A}$ is the right singular vector. In Eq. (34), $m$ denotes the number of modes, which is smaller than $l$ if the rank of $\boldsymbol{U}$ is smaller than $l$. The selection of $m$ is based on the criterion formulated in Eq. (30), that is, the first $m$ modes contain most of the energy of the whole displacement field. Therefore, the approximate displacement which can well characterize the complete displacement field is obtained according to $\Sigma$. By substituting the obtained displacement vector into Eq. (30), the approximated ESLs of dynamic problems is obtained:

$$\boldsymbol{F} = \{\boldsymbol{f}_1 \; \boldsymbol{f}_2 \ldots \boldsymbol{f}_m\} = \boldsymbol{K} \times \{\boldsymbol{\varphi}_1 \; \boldsymbol{\varphi}_2 \ldots \boldsymbol{\varphi}_m\} \tag{35}$$

### *4.2. Optimization problem formulation*

The general topology optimization problems based on static response in the density method are as follows:

| | | |
|---|---|---|
| Find | $\boldsymbol{b} \in R^n$ | (36a) |
| To minimize | $\boldsymbol{f}^T \boldsymbol{d}$ | (36b) |
| Subject to | $\boldsymbol{K}\boldsymbol{d} - \boldsymbol{f} = 0$ | (36c) |
| | $v^T \boldsymbol{b} \leq V$ | (36d) |
| | $0 < \boldsymbol{b}_{min} < \boldsymbol{b}_i < 1, i = 1, \ldots, n$ | (36e) |

where $\boldsymbol{b}$ is the artificial density of finite elements (design variable), $n$ the number of elements or design variable, $\boldsymbol{f}^T \boldsymbol{d}$ the compliance and $v$ the volume vector of the finite elements. $V$ is the volume fraction of the structure defined by the designer, $\boldsymbol{b}_i$ is the $i^{th}$ design variable and $\boldsymbol{b}_{min}$ is the lower bound of the design variables.

When dynamic topology response optimization is performed, the objective compliance and state equation constraints will change. The compliance is evaluated in the time domain, which is the summation of the compliance at all-time intervals. So, the dynamic topology optimization problem is described below:

$$\text{Find} \quad \boldsymbol{b} \in R^n \quad (37a)$$

$$\text{To minimize} \quad \sum_{u=1}^{l} c(t_u) \quad (37b)$$

$$\text{Subject to} \quad \boldsymbol{M}\ddot{\boldsymbol{d}}(t) + \boldsymbol{K}\boldsymbol{d}(t) - \boldsymbol{F}(t) = 0 \quad (37c)$$

$$v^T \boldsymbol{b} \leq V \quad (37d)$$

$$0 < \boldsymbol{b}_{min} < \boldsymbol{b}_i < 1, i = 1, \dots, n \quad (37e)$$

where $c$ represents the compliance at each time interval.

As mentioned earlier in Section 4.1, the exact ESLs is obtained as follows:

$$\boldsymbol{F}_{eq} = \boldsymbol{F}(t) - \boldsymbol{M}(\boldsymbol{b})\ddot{\boldsymbol{d}}(t) = \boldsymbol{K}\boldsymbol{d}(t) \quad (38)$$

Thus, based on the ESL method, the dynamic response topology optimization is performed by dynamic simulation and static topology optimization under multiple static loads. And topology optimization problem in Eq. (37) can be transformed multiple static topology optimizations as follows:

$$\text{Find} \quad \boldsymbol{b} \in R^n \quad (39a)$$

$$\text{To minimize} \quad \sum_{u=1}^{l} \left(\boldsymbol{F}_{eq}(t_u)\right)^T \boldsymbol{d}(t_u) \quad (39b)$$

$$\text{Subject to} \quad \boldsymbol{K}\boldsymbol{d}(t_1) - \boldsymbol{F}_{eq}(t_1) = 0$$

$$\boldsymbol{K}\boldsymbol{d}(t_2) - \boldsymbol{F}_{eq}(t_2) = 0$$

$$\dots \quad (39c)$$

$$\boldsymbol{K}\boldsymbol{d}(t_l) - \boldsymbol{F}_{eq}(t_l) = 0$$

$$v^T \boldsymbol{b} \leq V \quad (39d)$$

$$0 < \boldsymbol{b}_{min} < \boldsymbol{b}_i < 1, i = 1, \dots, n \quad (39e)$$

It can be seen from Eq. (39c) that the calculation cost of static topology optimization is positively proportional to the number of time intervals. To decrease the computational cost of static topology optimization, according to the reduction method for multiple ESLs in Section 4.1, the objective function can be given as:

$$\sum_{u=1}^{m}(K\boldsymbol{\Phi}_u)^T \boldsymbol{d}_u \tag{40}$$

Constraint in (39c) becomes:

$$\boldsymbol{K}\boldsymbol{d}_1 - \boldsymbol{f}_1 = 0$$
$$\boldsymbol{K}\boldsymbol{d}_2 - \boldsymbol{f}_2 = 0$$
$$\ldots$$
$$\boldsymbol{K}\boldsymbol{d}_m - \boldsymbol{f}_m = 0 \tag{41}$$

where $m$ denotes the number of reduced ESLs, that is, approximated ESLs, $\{\boldsymbol{d}_1\ \boldsymbol{d}_2\ \ldots\ \boldsymbol{d}_m\}$ the displacement vector corresponding to $\{\boldsymbol{f}_1\ \boldsymbol{f}_2\ \ldots\ \boldsymbol{f}_m\}$. From the examples studied later, $m$ is much less than $l$, that is, $m \ll l$, the efficiency of static topology optimization is significantly improved.

### 4.3. Sensitivity analysis & Variable update

Based on the FEM, the structure is discretized into standard square elements. In this study, the stiffness matrix is constructed by the stiffness interpolation function based on SIMP method [34, 35], i.e. each element $e$ is assigned a density $\boldsymbol{b}_e$, so the interpolation formula of the material is obtained:

$$\boldsymbol{E}_e(\boldsymbol{b}_e) = \boldsymbol{E}_{min} + \boldsymbol{b}_e^p(\boldsymbol{E}_0 - \boldsymbol{E}_{min}), \quad \boldsymbol{b}_e \in [0,1] \tag{42}$$

where $\boldsymbol{E}_0$ is the element stiffness of the material, $\boldsymbol{E}_{min}$ is a very small stiffness assigned to void element in order to avoid the singularity of the total stiffness matrix of the assembly, and $p$ is a penalization factor (typically $p = 3$) introduced to ensure a clearer structural topology.

In order to obtain clear topology optimization result, the volume interpolation function based on the Heaviside projection function is used to construct the mass matrix.

$$V_e(\boldsymbol{b}_e) = 1 - e^{-\boldsymbol{b}_e p\_0} + \boldsymbol{b}_e e^{-p\_0} \tag{43}$$

where $p\_0$ is a penalization factor for Heaviside projection.

Thus, the interpolation formula of the material based on the Heaviside projection is obtained:

$$\boldsymbol{E}_e(\boldsymbol{b}_e) = \boldsymbol{E}_{min} + V_e^{\,p}(\boldsymbol{E}_0 - \boldsymbol{E}_{min}), \quad \boldsymbol{b}_e \in [0,1] \tag{44}$$

In this way, the objective function of topology optimization problem can be rewritten in the following form:

$$z(\boldsymbol{b}) = \sum_{u=1}^{m}\left(\sum_{e=1}^{n} \boldsymbol{E}_e(\boldsymbol{b}_e)\boldsymbol{d}_e^{\,T}\boldsymbol{k}_0\boldsymbol{d}_e\right)_u \tag{45}$$

where $z$ is the compliance, $\boldsymbol{d}_e$ is the element displacement vector, $k_0$ is the element stiffness matrix of unit Young's modulus.

The sensitivities of the objective function $z$ and the material volume $V$ with respect to the element densities are given by:

$$\frac{\partial z}{\partial \boldsymbol{b}_e} = \sum_{u=1}^{m}(-pV_e^{\,p-1}(e^{-\boldsymbol{b}_e p0}p0 + e^{-p0})(\boldsymbol{E}_0 - \boldsymbol{E}_{min})\boldsymbol{d}_e^{\,T}\boldsymbol{k}_0\boldsymbol{d}_e)_u \tag{46}$$

$$\frac{\partial V}{\partial \boldsymbol{B}_E} = (1 - \boldsymbol{E}_{min})(e^{-\boldsymbol{b}_e p0}p0 + e^{-p0}) \tag{47}$$

The optimization problem in the section 4.2 is solved by means of a standard Optimality Criteria (OC) method. The update and transfer of design variables are shown in the optimization loop, as shown in Algorithm 1, that is, the overall process of the proposed dynamic topology optimization strategy.

**Algorithm 1:** Dynamic topology optimization

**Input**: $b^0, Tol_{dy}, MaxIter$
**Output**: Design variables $b$

| | |
|---|---|
| 1 | Initialization: $k=1, b^k = b^0$ |
| 2 | **while** $Max(b^{k+1} - b^k) > Tol_{dy}$ && $k < MaxIter$ **do** |
| 3 | Calculate the approximate displacement response $U$ of dynamic problems by the OSDCA as in Eqs. (14-26) |
| 4 | Calculate the principal displacement $\Phi = POD(U) = \{\varphi_1\ \varphi_2\ ...\ \varphi_m\}$ as in Eq. (33) |
| 5 | Calculate ESL sets $F_{eq}$ as in Eq. (35) |
| 6 | Execute the multi-static topology optimization process with Eq. (40) as the objective and Eq. (41) as the state equation constraint |
| 7 | Update design: $k = k+1$ |
| 8 | **End** |

## 5. Numerical case studies

Several numerical examples were analysed to evaluate the proposed method in this section. To get insight into the performance of optimization framework, 2D and 3D structures under different load excitations and different model sizes are investigated. Additionally, the influence of some important parameters on the optimization effectiveness and results is illustrated by comparing the Model Superposition Method (MSM) with the proposed method. The material is set to steel with elastic modulus $E$=201GPa, density $\rho$ =7.85g/cm$^3$, and the Poission's ratio $v$=0.33. For the transient analysis, a Newmark-$\beta$ scheme was chosen to solve the time-dependent equation.

### 5.1. Cantilever beam with a hole

This example mainly illustrates the performance of the dynamic reanalysis method in terms of computational efficiency and accuracy. Through comparative study with MSM, its acceleration performance is verified. In order to further reflect the results of dynamic topology optimization, the corresponding static topology optimization problems are also compared. A 2D cantilever structure with a hole is investigated in this example, whose thickness is set as 0.01$m$. The length is 2 $m$ and the width is 1 m ($L$ = 1). The left end is fixed and two dynamic forces are enforced on the upper right corner and lower middle of the structure respectively, as shown in the Figure. 4 (a). The external load curve is shown in Figure. 4 (b).

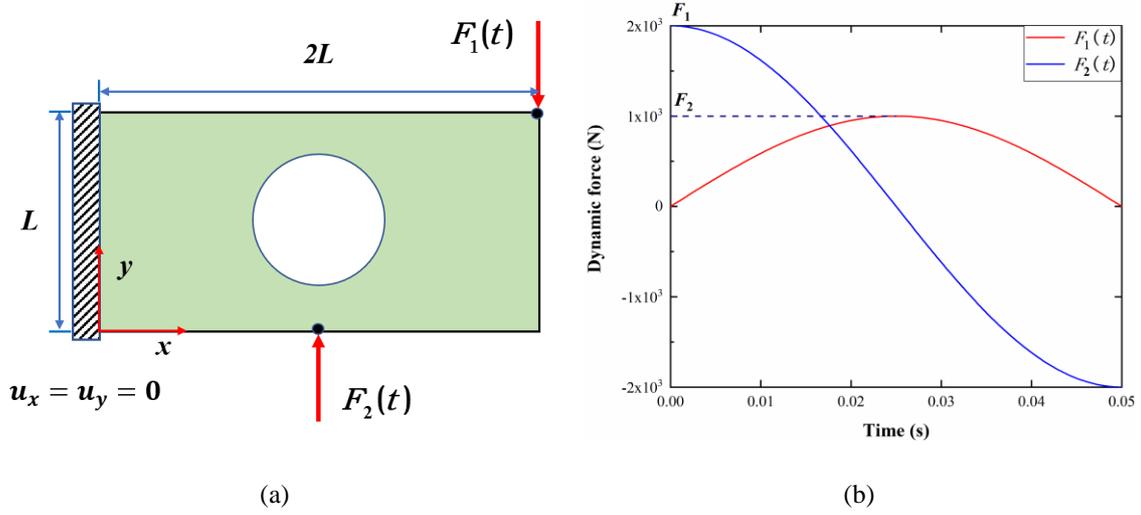

(a)            (b)

Fig. 4 The cantilever beam: (a) domain and boundary conditions and (b) applied load

In order to obtain a clear and convergent topology optimization result, the strategy of combination of the dynamic adjustment of penalty factors and projection technique is used. The specific setting is: the smaller penalization factor is given at the early stage of optimization, in which the density design variables are updated sharply; as topology optimization progresses, the update of design variables gradually stabilizes, and the given penalization factor increases, as shown in Fig. 5.

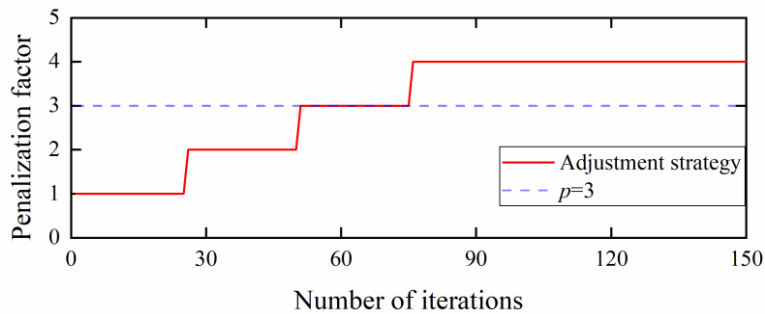

Fig. 5 The adjustment strategy of penalization factor

Additionally, the Heaviside projection method defined in Eq. (43) is applied throughout the topology iteration process to obtain a clearer topology configuration. Table 1 shows the comparison of optimization results before and after applying the parameter dynamic adjustment strategy. From the comparison results, it can be found that the result obtained by the suggested strategy yields more clear results.

Table 1 The evolution of structure topology with different optimization strategies

| Iterative step | Penalization factor ($p = 3$) | Adjustment strategy with Heaviside projection |
|---|---|---|
| 10 | | |
| 20 | | |
| 50 | | |
| 150 | | |

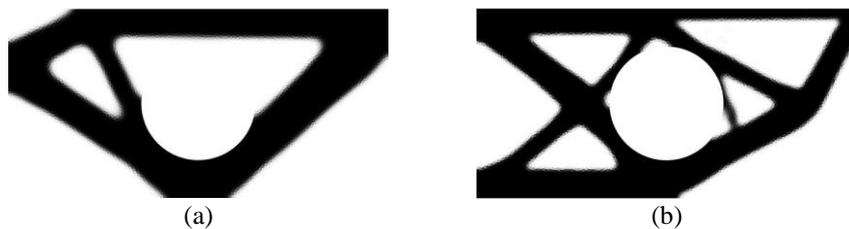

Besides, for the sake of reflecting the dynamic characteristics of the studied problem, the static and dynamic topology optimization are performed. The static topology optimization problem with the same design domain and boundary is presented in Fig. 4 (a), and the dynamic load is replaced by the static load ($F_1$ and $F_2$). Fig. 6 (a) presents the result of static topology optimization and Fig. 6 (b) is the result of dynamic topology optimization. It can be found that static and dynamic problems show significantly different topology optimization configurations.

(a)          (b)

Fig. 6 The comparison of topology optimization: (a) static problem; (b) dynamic problem.

(c) the evolution of first order fundamental natural frequency for dynamic structure

After determination of the dynamic topology optimization problem, the performance of the proposed method in optimization precision and calculation efficiency is mainly discussed.

In each dynamic simulation, dynamic reanalysis is used to improve computational efficiency, and the activation of reconstructing basis vectors is dynamically adjusted based on the changes in adjacent dynamic topology iteration structures. The variation of structural stiffness matrix ($\Delta K$) is defined to represent the structural change of dynamic topology optimization iteration, as follow:

$$\Delta K = \left|\frac{K^{k+1} - K^k}{K^k}\right| \tag{48}$$

When $Max(\Delta K)$ is greater than a given threshold ($Tol_{rb}$), the basis vectors in dynamic reanalysis are reconstructed to ensure the precision of the approximate response. Besides, the number of basis vectors in the reduced space should be determined for the reanalysis. To this end, the performance of topology optimization with different number of basis vectors are tested. Table 2 lists the iterative topology optimization results of approximate models with different number of basis vectors. Compared with the results of full model in Table 1, the approximate model constructed based on 3 basis vectors can obtain a basically consistent optimal configuration.

Table 2 The topology optimization results with different number of basis vectors

| Number of basis vectors (s) | 10th iteration | 20th iteration | 50th iteration | 150th iteration |
|---|---|---|---|---|
| s=1 | 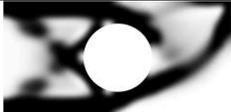 | 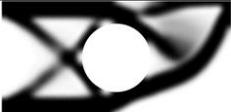 | 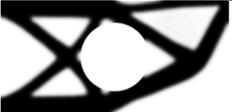 | 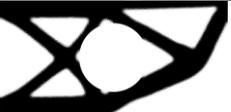 |
| s=2 | 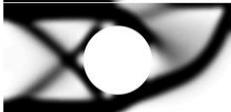 | 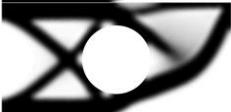 | 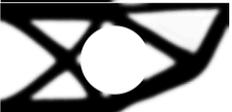 | 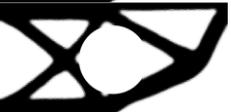 |
| s=3 | 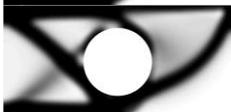 | 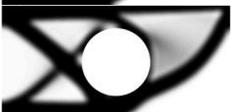 | 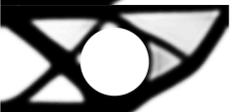 | 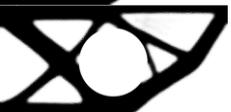 |

Further, since not every two adjacent steps in the entire dynamic topology optimization iteration process changes a lot, the strategy of adaptive reconstruction of reduced space (basis vectors) is performed. Fig. 7 shows the structural change histories of approximate models in dynamic topology optimization. At the early stage of optimization, the reduced space is reconstructed at each iteration due to significant structural changes. The red dotted line is the given threshold ($Tol_{rb}=0.01$). When the structural change is lower than the threshold, it is not necessary to reconstruct the reduced space in pace with optimization. The bar chart explains the implementation of adaptive strategies for different approximation models. It can be found that the approximate model with 1 or 2 basis vectors needs more reconstruction of reduced space in the optimization process, and the structure of adjacent iterations at the early stage of optimization also changes greatly. Instead, when the number of basis vectors is 3 or more, the model approximate precision is improved, thus the number of the basis vectors reconstruction is reduced, achieving more efficient calculation.

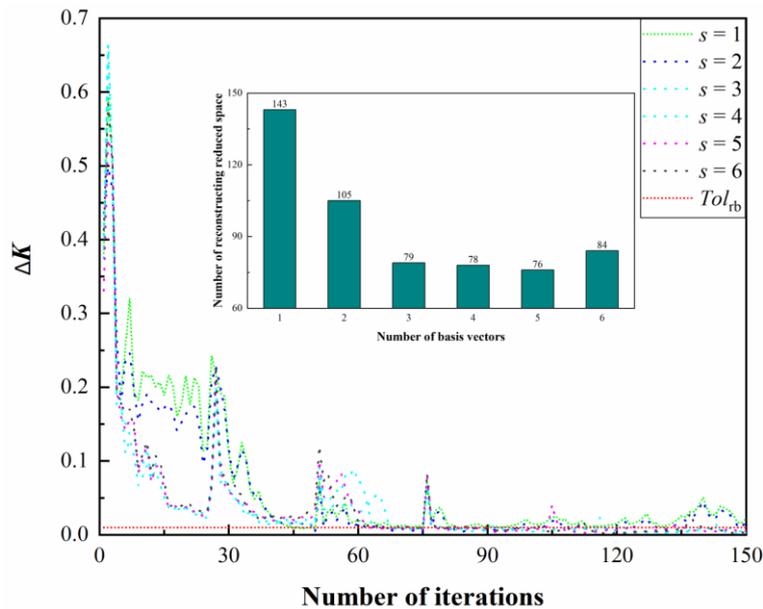

Fig. 7 The evolution of structural change for different models

Considering the above results, further comparative research will be conducted on the MSM and the proposed reanalysis method. For calculation accuracy, 3 basis vectors are given for reanalysis and 3 modes for MSM, besides, the displacement in the *y*-direction at the loading point of dynamic force $F_1$ is compared with complete analysis (Newmark) respectively.

To sufficiently certificate the effectiveness of the approximate method, the initial and modified structure are designated as the adjacent iterations with the maximum structural change according to Fig. 7. Figure 8 shows the approximate displacement of the modified structure and the error compared to the complete analysis. Apparently, it can be known that the proposed method outclasses MSM in calculation accuracy, in which the displacement error is maintained on $10^{-12}$ order of magnitude. On the contrary, the approximate displacement solved by MSM differentiates significantly from the exact solution in the first 0.01 seconds, which has a maximum error of nearly 100%.

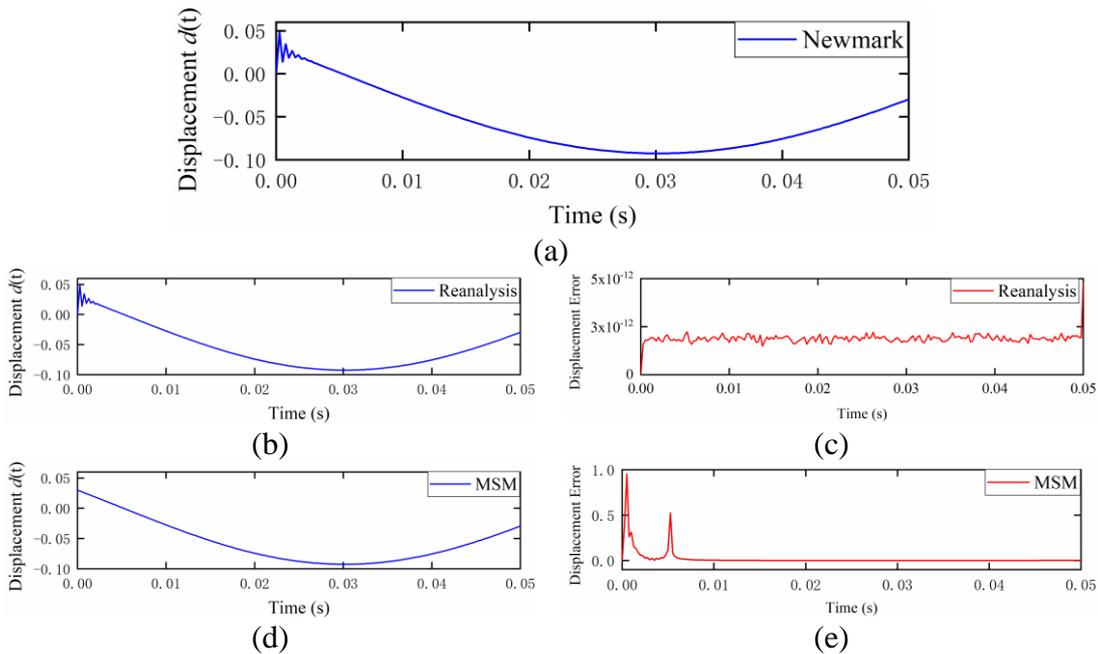

Fig. 8 The approximate displacement and error for OSDCA and MSM:
(a) the displacement of complete analysis; (b) the approximate displacement for OSDCA;
(c) the error for OSDCA; (d) the approximate displacement for MSM; (e) the error for MSM.

For computational efficiency, the cantilever beam model is divided into 10000 elements, containing 19683 nodes, and the total degree of freedom is 39366. Besides, the dynamic simulation time is 0.05 seconds and divided into 200 time intervals. The acceleration effect of OSDCA and MSM is comparative studied on account of the above finite element model, that each dynamic simulation time of the entire optimization evolution is tested, as described in Fig. 9. From the curves, it can be seen that dynamic simulations take sharply less time with OSDCA than MSM, thus it has a significant acceleration effect.

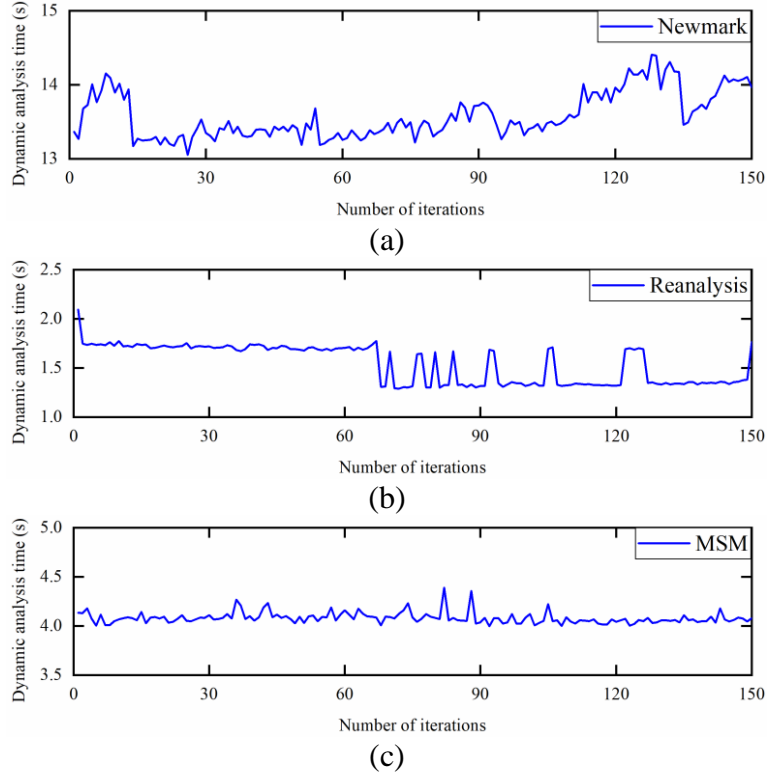

Fig. 9 Each dynamic analysis time for optimization evolution with different methods:
(a) Newmark; (b) OSDCA; (c) MSM.

In view of the above results, the OSDCA implemented model reduction in spatial domain for dynamic simulations. Then, the time domain reduction is performed by the displacement's POD approximation. Given that the classical ESL method, the approximated ESLs are calculated instead of exact ESLs. Further, stemmed from the ratio definition in Eq. (30), the number of approximated ESLs ($m$) is determined by the singular value of dynamic displacement snapshots, as follow:

$$\varepsilon(m) = \frac{\sum_{j=1}^{m} S_j}{\sum_{j=1}^{l} S_j} \qquad (30)$$

where $S$ is the singular value, $S = [S_1, S_2, ..., S_j]$, which are sorted.

Apparently, the dynamic displacement field should be changing owing to the constantly changing structure in the dynamic topology evolution. Thus, by setting a specific energy ratio (here $\varepsilon = 0.9$), the approximated ESLs of each dynamic optimization iteration are adaptively constructed to achieve model reduction in the time domain.

Finally, the dynamic topology optimization results of approximated model, which is reduced in both space and time domains, are discussed. Considering the dynamic problem studied, the exact number of ESLs is 200. Fig. 10 (a) shows the number of approximated ESLs adaptively determined based on the given energy ratio, which is visibly smaller than 200. Thus, the static topology optimization should be obviously reduced. Fig. 10 (b) and (c) are the static topology optimization time with exact and approximated ESLs respectively. Also, Fig. 10 (c) displays the generation time of approximated ESLs. Although an additional operation, that calculate the dynamic displacement's POD approximation, is involved for static topology optimization with approximated ESLs, the combined time of generation and optimization illustrated in Fig. 10 (c) is much smaller than exact ESLs. Thus, the efficiency of static topology optimization is significantly improved.

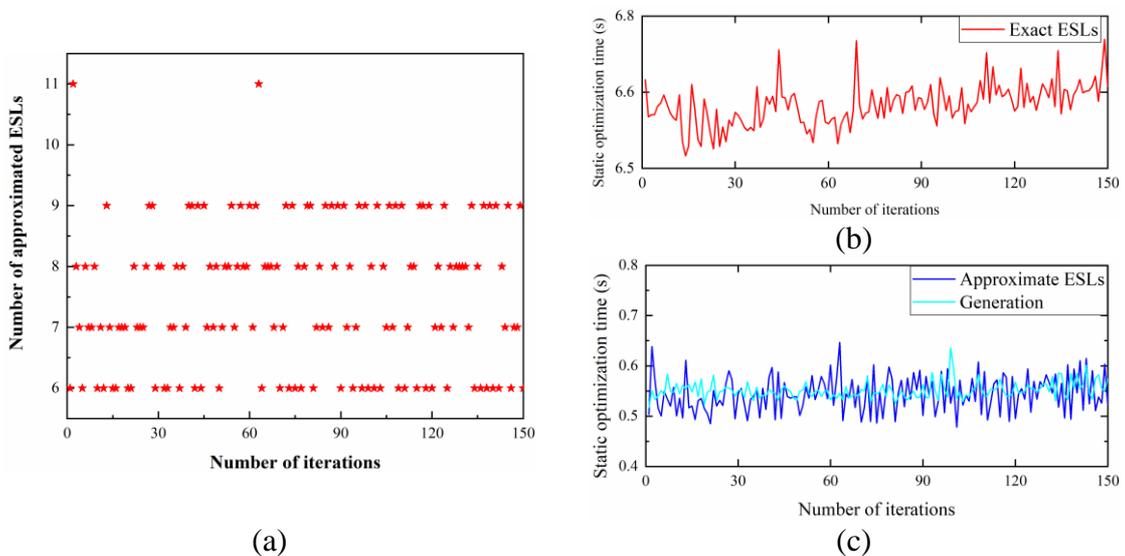

Fig. 10 Results of the approximated ESLs: (a) number; (b) static topology optimization time with exact

ESLs; (c) time of generation and static topology optimization for approximated ESLs.

In addition to the efficiency improvement, the performance of the optimized structure is verified. Fig. 11 (a) shows the displacement in the y-direction at the loading point of dynamic force $F_1$ of the initial and the optimal structure. It can be found that the deformation of the optimal structure is smaller than the initial structure in the entire time history, it indicates the improvement of the structural stiffness. Furthermore, the iterative curves of first order fundamental natural frequency and some procedures of topological configuration are illustrated in Fig. 10 (b), from which it is remarkable that: compared with the initial structure, the fundamental natural frequency increased for optimized structure; utilizing the reduced strategies for dynamic topology optimization, the optimal configuration is coincident with full model.

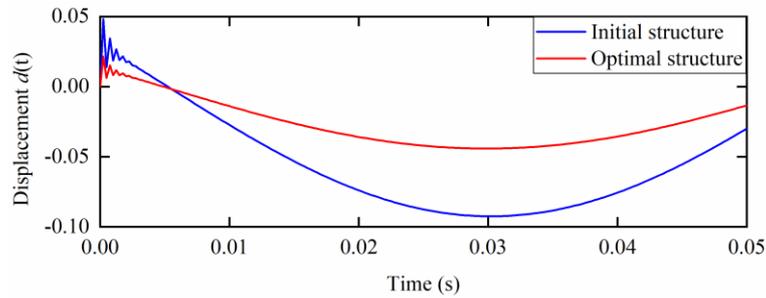

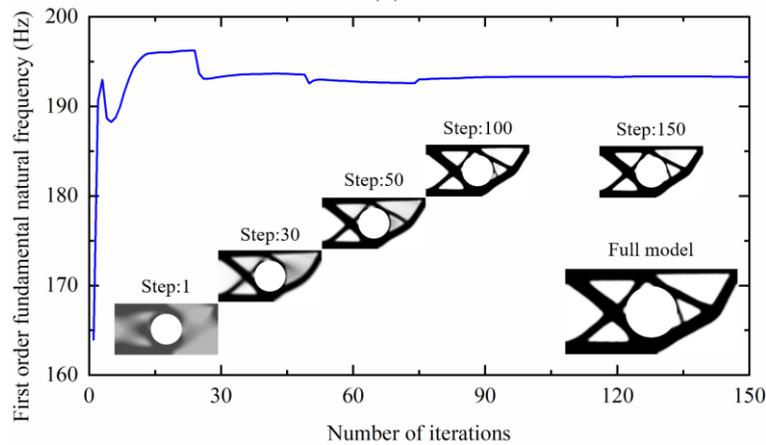

Fig. 11 The results of dynamic topology optimization: (a) displacement; (b) first order fundamental natural frequency and optimal configurations.

## 5.2. Bridge with moving load

This example aims to verify the acceleration performance of the presented methodology for two-dimensional structure. To this end, a bridge structure is further investigated as shown in Fig. 12. The bridge with length $L = 30$ m and height $H = 13.5$ m is subjected to a uniform load $F = 100$ kN of length $L_0 = 5$ m that moves over the deck of the bridge at a constant speed $v_0 = 30$ km/h.

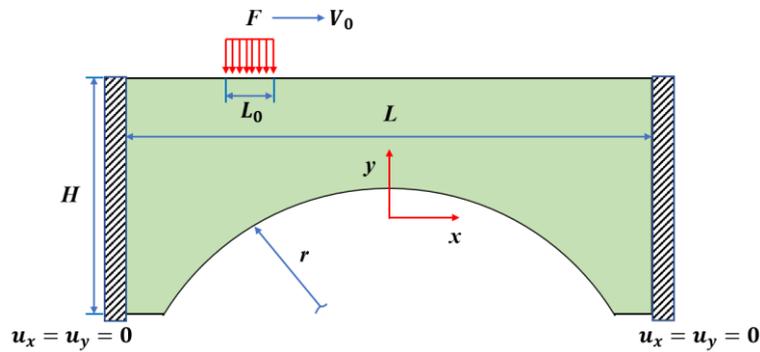

Fig. 12 The bridge structure and boundary condition

The number of basis vectors should be firstly determined for construction of a reduced space for bridge design problems. Fig. 13 illustrates the optimal bridge design results, which is clearly concluded that the reduced model with 6 basis vectors can obtain consistent design with full model.

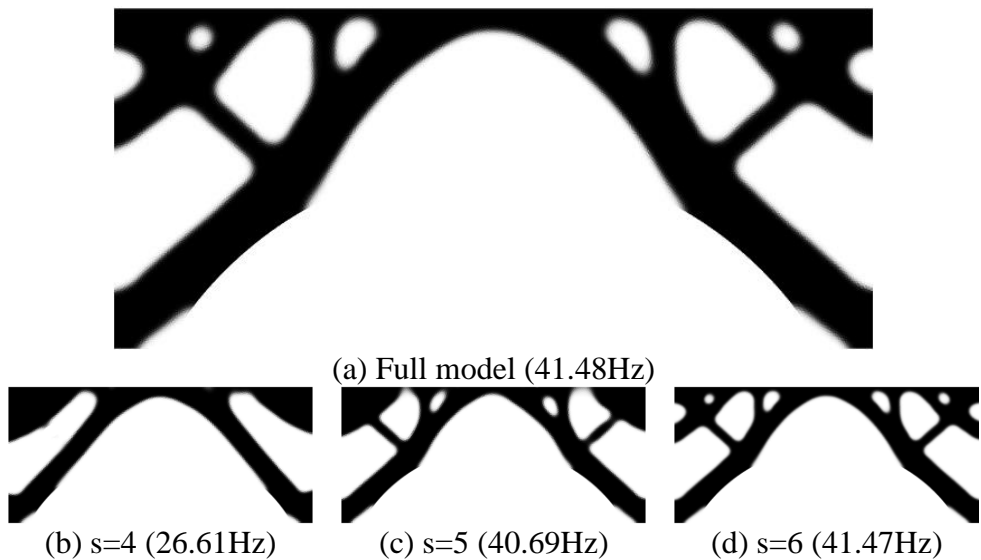

(a) Full model (41.48Hz)

(b) s=4 (26.61Hz)　　(c) s=5 (40.69Hz)　　(d) s=6 (41.47Hz)
Fig. 13 The topology optimization results of different reduced model

Then, in view of the determined reduced space dimension, bridge structures with variant scales are optimized, and the information of their finite element models is listed in Table 3.

Table 3 The finite element model information

| Bridge model | 1 | 2 | 3 | 4 | 5 |
|---|---|---|---|---|---|
| Node | 19704 | 39372 | 58965 | 78639 | 98320 |
| Element | 10000 | 20000 | 30000 | 40000 | 50000 |
| DOF | 39408 | 78744 | 117930 | 157278 | 196640 |

In each dynamic topology optimization iteration, dynamic simulation is firstly performed to attain the displacement. Then, its POD approximation is calculated for constructing approximated ESLs. Finally, the bridge structure subjected to these multiple ESLs is optimized. To this end, the acceleration performance of the proposed methodology is proofed in these three solution segments respectively.

(1) Dynamic analysis.

Fig. 14 (a) shows the dynamic calculation time of full model, given in Table 3, for each optimization iteration. Correspondingly, the reduced model constructed by OSDCA are illustrated in Fig. 14 (b). It can be concluded that: implementing the OSDCA, the dynamic simulation time is significantly declined; besides, the history curves of each dynamic simulation time over the entire optimization verify the adaptive reconstruction reduced basis vectors strategy. Notably, the time is sharply reduced where there is no need for reconstruction.

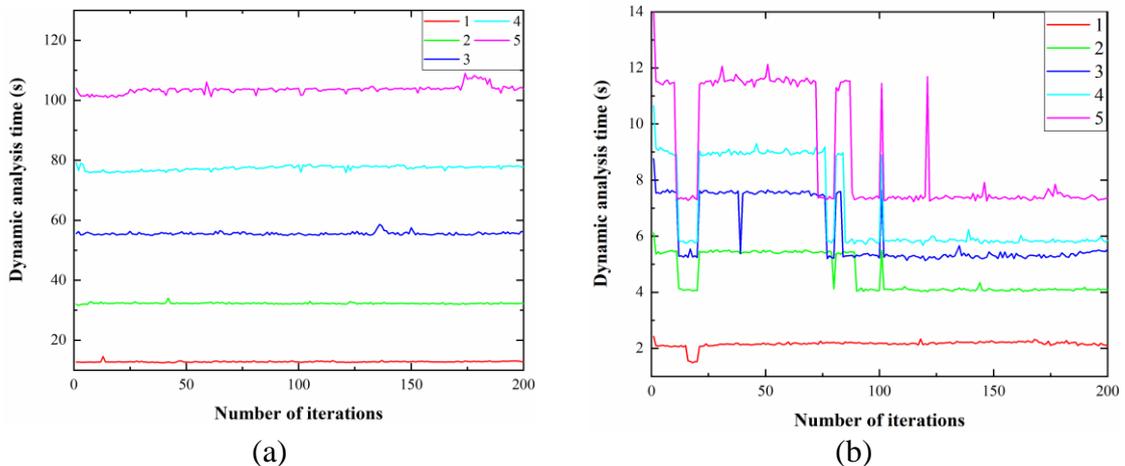

(a)            (b)

Fig. 14 The dynamic analysis time history in optimization for different structure scales:

(a) Newmark; (c) OSDCA.

Owing to the introduction of adaptive strategy, the dynamic analysis time in the whole dynamic topology optimization process is considered. So, considering the results revealed in Fig. 14, each dynamic simulation time should be summed and the speedup ratio is calculated for reduced models with different scales, as shown in Fig. 15. It is remarkable that as the model scale increases, the speedup ratio tends to increase, indicating that this method has more obvious potential for large-scale problems.

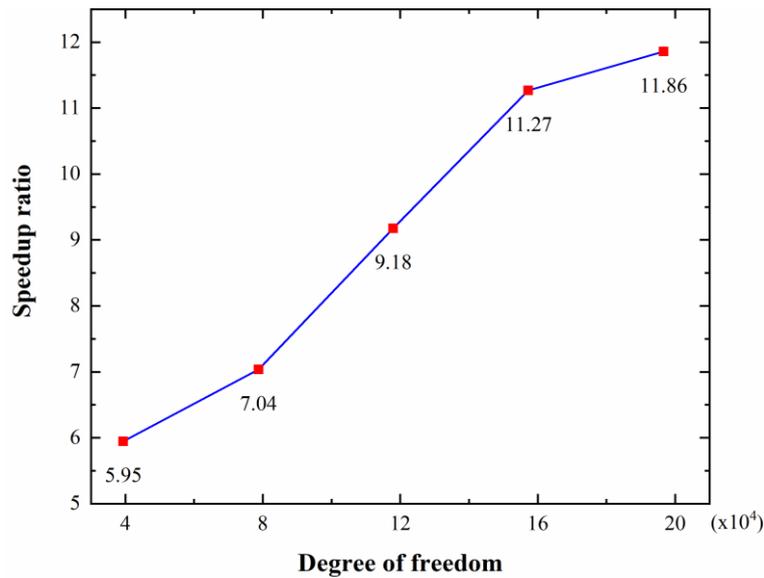

Fig. 15 The speedup ratio for OSDCA

(2) Generation of the approximated ESLs.

This portion connects dynamic solution and static topology optimization, involving the matrix decomposition of displacement field, which adds extra calculation time compared to full model. Furthermore, with the increase of model scale, its time-consuming increases linearly, displayed in Fig. 16 (a). Nevertheless, its time proportion is too marginal to hinder the application to the large-scale problems, which can be clearly illustrated in Fig. 16 (b).

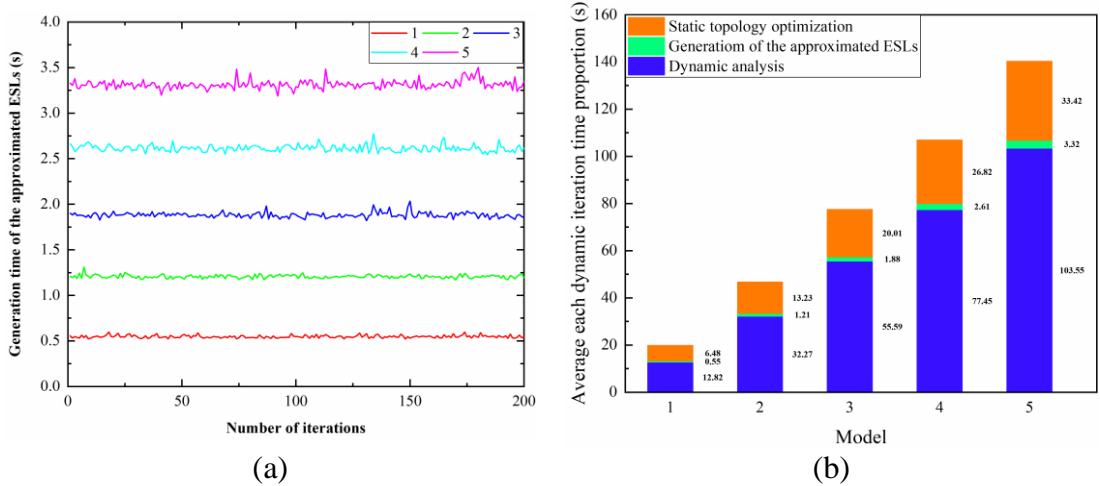

(a)                          (b)

Fig. 16 The approximated ESLs: (a) generation time; (b) time proportion for full model.

(3) Bridge design optimization with static loads.

Fig. 17 displays the time of bridge optimization with exact and approximated ESLs respectively. Owing to adaptive construction of the approximated ESLs based on energy ratio ($\varepsilon = 0.9$), the number of ESLs contained in each iteration is constantly changing, which is reflected in optimization time shown in Fig. 17 (b). Besides, generally speaking, the number of approximate ESLs is smaller than that of exact ESLs, which is consistent with total time intervals, so that the capacity is promoted. Similarly, the average each optimization time is computed, as displayed in Fig. 17 (c). It is revealed that as the model scale increases, the acceleration affect is maintained, indicating that this method accommodates to large-scale problem.

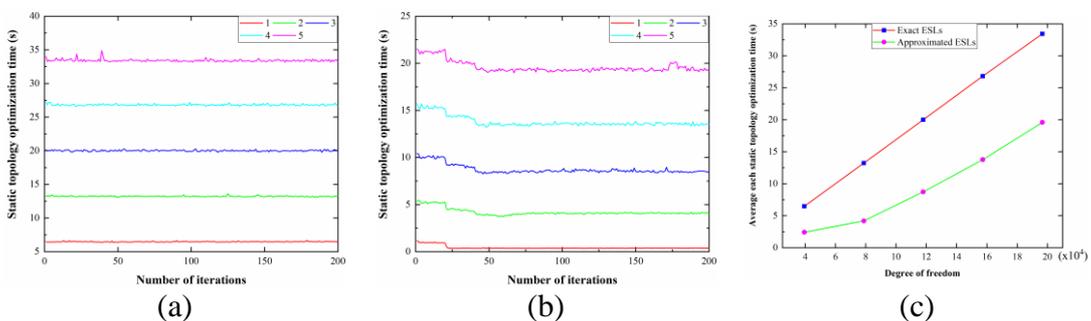

(a)                  (b)                  (c)

Fig. 17 The static topology optimization time history for different structure scales:
(a) exact ESLs; (b) approximated ESLs; (c) average each time.

Finally, the time proportion for reduced model is represented in Fig. 18. In each dynamic iteration, although total solution time is significantly reduced compared to full model illustrated in Fig. 16 (b), the acceleration effect for static topology optimization is unconspicuous, further, its time proportion would magnify as structural scale increases, pulling the overall acceleration capacity down.

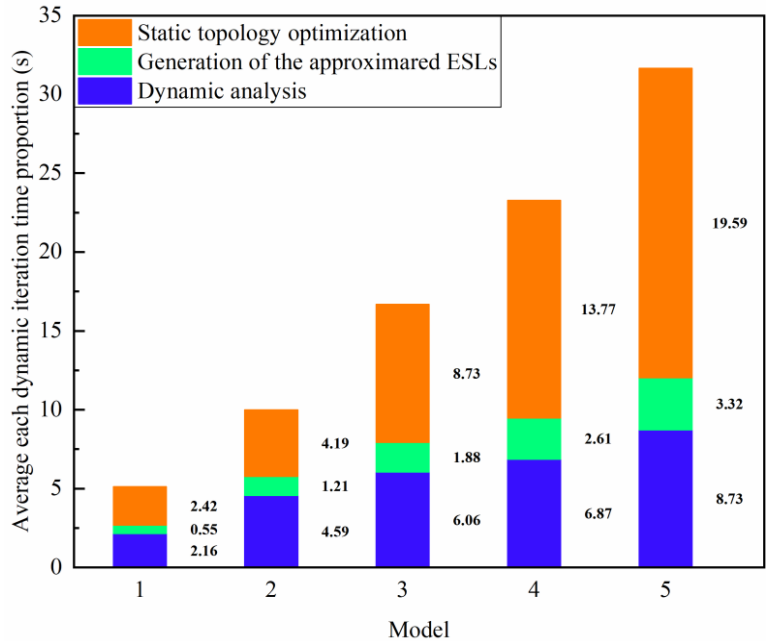

Fig. 18 The time proportion for reduced model.

Given the above conclusion, it is necessary for exploring the parameter $\varepsilon$, which directly impacts on static optimization time and optimal configuration. Thus, different energy ratios are given on Model. 3 calculated with 200 time intervals. The results are illustrated in Fig. 19, which can be visibly concluded that on the premise of pledging the optimal configuration, appropriately reducing the energy ratio can significantly reduce the number of approximate ESLs, thereby improving computational efficiency.

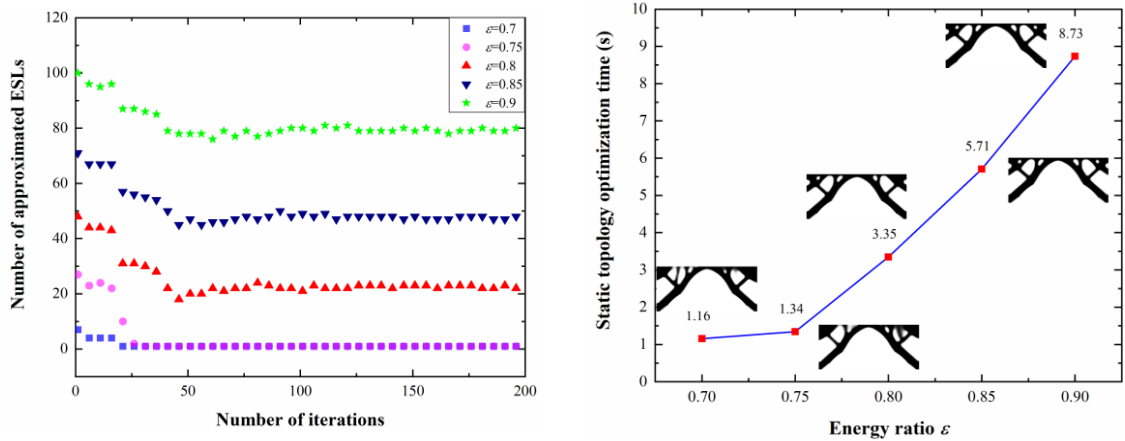

(a)                              (b)

Fig. 19 The results of different energy ratios:

(a) number of approximated ESLs; (b) time and topology configuration.

## 5.3. 3D structure under random load

The 3D structure is further investigated to demonstrate the validity and feasibility of the developed reduced strategy. To this end, a hook structure with high $H = 0.15$ m, width $W = 0.1$ m and thickness 0.06 m is subjected to randomized dynamic force, as displayed in Fig. 20. The upper half hole wall surface of the hook structure is constrained, and randomized dynamic loads are applied to the lower hook surface.

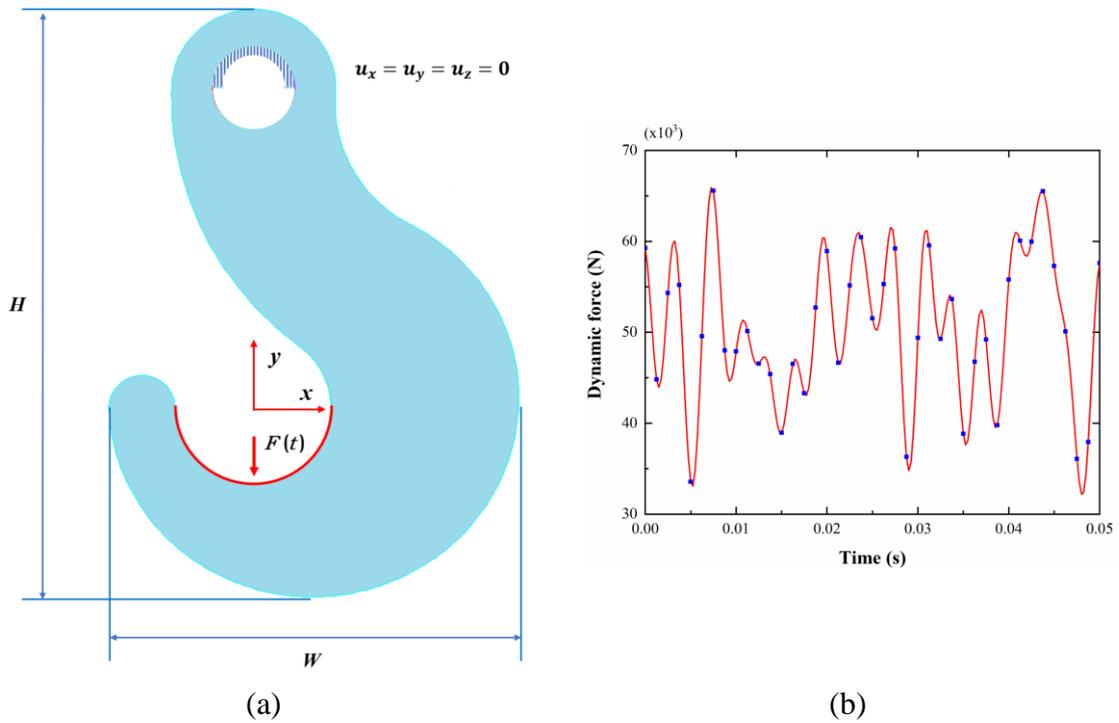

(a)                              (b)

Fig. 20 The hook structure: (a) domain and boundary conditions and (b) applied dynamic load

Firstly, the OSDCA is performed to calculate the dynamic displacement, in which a given threshold ($Tol_{rb} = 0.01$) is presumed. Besides, the reduced basis vectors are determined by testing the reduced space with different number of basis vectors, as shown in Fig. 21. The reduced models are verified through the optimal configuration, and structural variation history further reflects the execution of adaptive strategy throughout dynamic topology optimization. Remarkably, the reduced model with 3 basis vectors shows optimal configuration with precision. In addition, the structural variation tends to stabilize and reaches the given threshold along with optimizing, which preferably implements the present OSDCA methodology to improve efficiency.

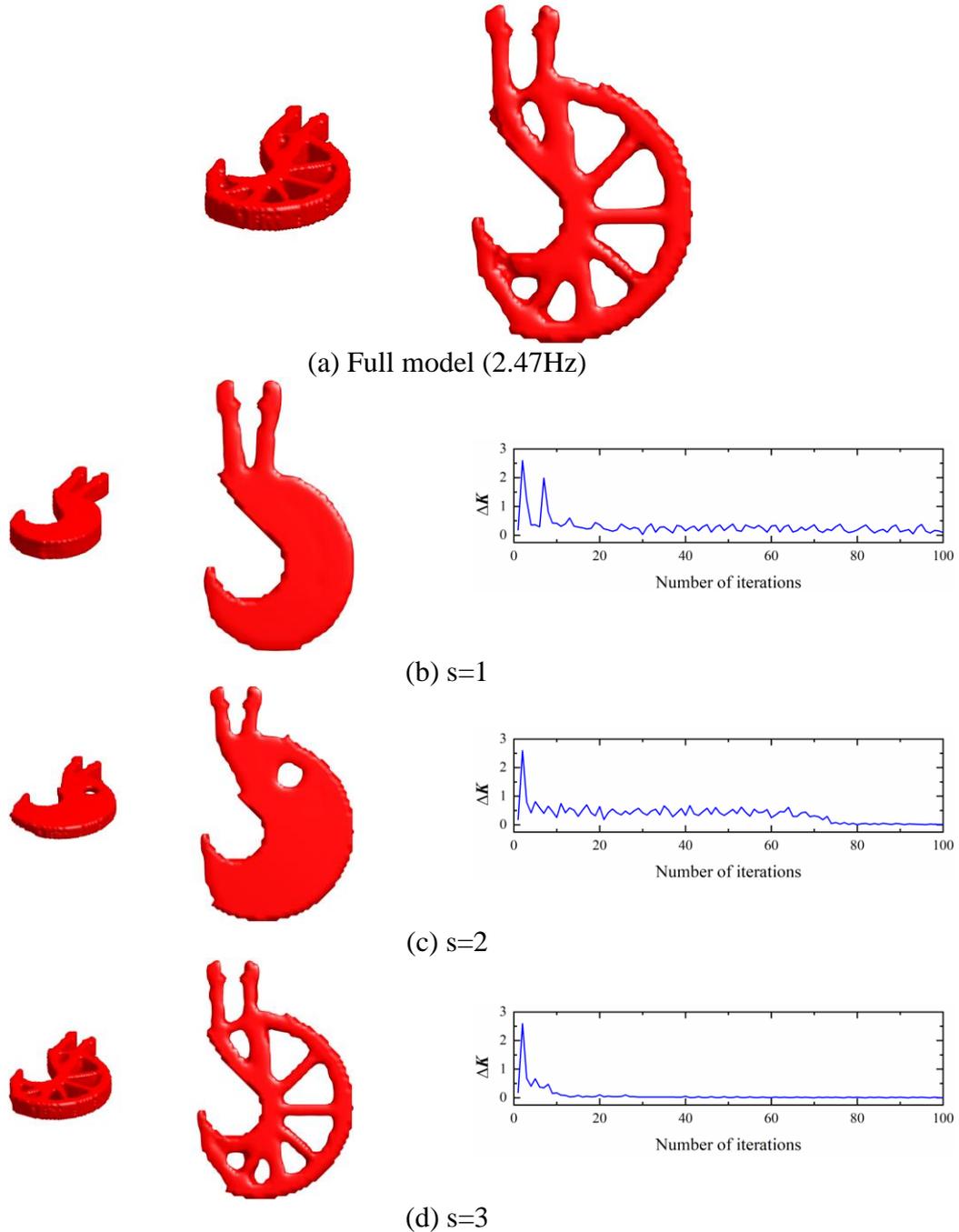

(a) Full model (2.47Hz)

(b) s=1

(c) s=2

(d) s=3

Fig. 21 The topology optimization results of different reduced model:

optimal configuration (left) and structure change history (right).

Then, the approximated ESLs are adaptively calculated in view of a specific energy ratio ($\varepsilon = 0.95$). Finally, the dynamic topology optimization results of approximated model for 3D structure are discussed, including performance and efficiency. Fig. 22 (a) shows the structural compliance evolution in the dynamic topology optimization. And it indicates that the performance of structural stiffness is gradually improved until stable, which is further validated that the deformation of where subjected to dynamic force decreases after optimization, as displayed in Fig. 22 (b). Additionally, this example is divided into 200 time intervals, thus, the exact ESLs is given.

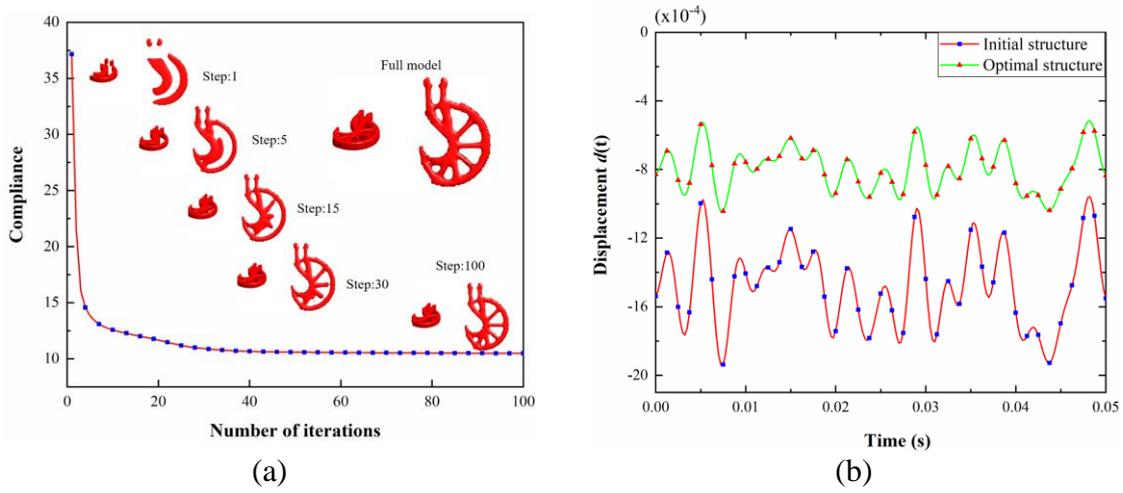

(a)                        (b)

Fig. 22 The results of dynamic topology optimization:
(a) the compliance evolution and topology configuration; (b) displacement.

Additionally, 3D hook structures with variant scales are optimized for exploring the acceleration performance of presented methodology, and their finite element model information are listed in Table 4.

Table 4 The finite element model information

| Hook model | 1 | 2 | 3 | 4 |
|---|---|---|---|---|
| Node | 11727 | 20774 | 32433 | 45973 |
| Element | 2508 | 5236 | 9450 | 15366 |
| DOF | 35181 | 62322 | 97299 | 137919 |

Given the discussion in Section 5.2, the acceleration performance of the presented methodology should be similarly proofed in two solution segments respectively for 3D problem.

According to results of dynamic analysis time demonstrated in Fig. 23, it is can be concluded that: (i) the calculation efficiency is enormously promoted; (ii) the adaptive strategy is triggered in the mid and later optimization stages, further improving its computational efficiency; (iii) with structural scale increasing, the acceleration capacity grows powerful, which is manifested in steadily increasing speedup ratio; (iv) further, the 3D dynamic analysis problem is better than 2D problem studied in Section 5.2 with regard to the acceleration capacity of presented methodology.

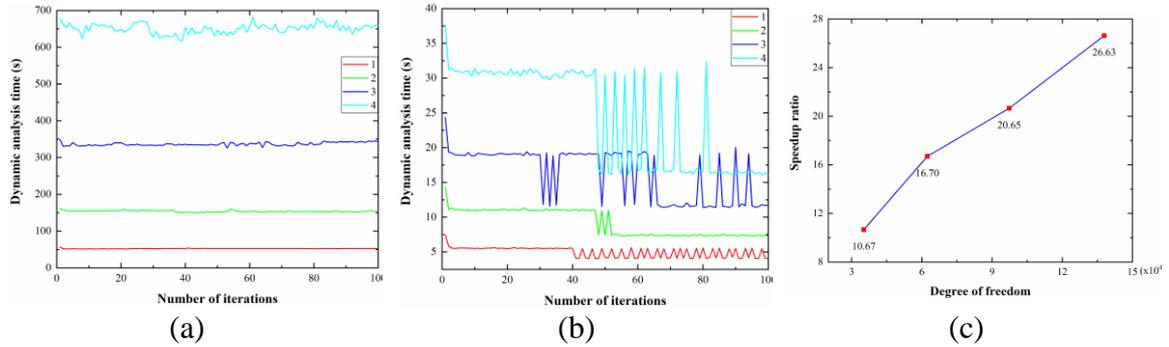

Fig. 23 The dynamic analysis time history in optimization for different structure scales:
(a) Newmark; (b) OSDCA; (c) speedup ratio.

Then, the hook structure optimization subjected to multiple ESLs is performed. And the average each optimization time is calculated, as displayed in Fig. 24. It is suggested that the acceleration capacity is maintained with structural scale increasing, which is stronger than 2D structure optimization ulteriorly.

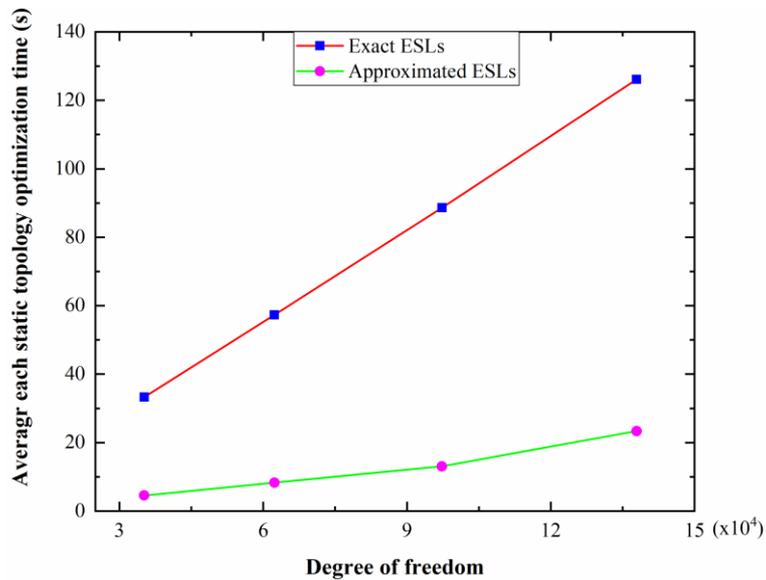

Fig. 24 The average each static topology optimization time ($\varepsilon = 0.95$)

Eventually, based on the presented methodology, the calculation time is significantly reduced. Furthermore, there is stronger acceleration potentiality on large-scale dynamic problem, which can be illustrated in Fig. 25.

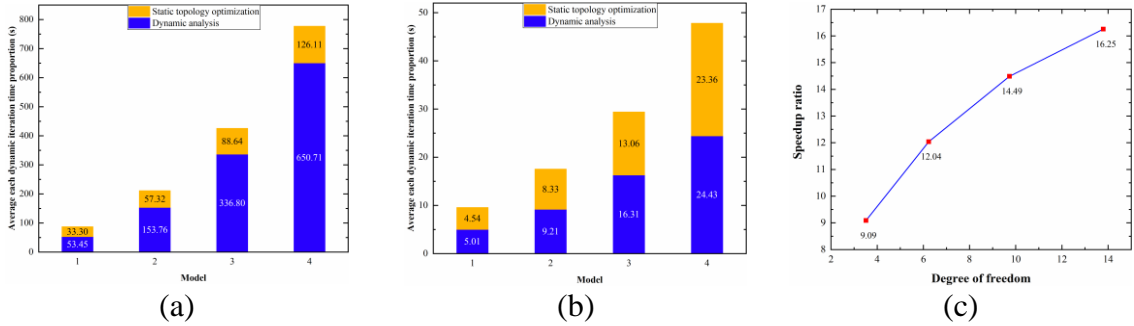

Fig. 25 The time proportion for average each dynamic iteration:
(a) full model; (b) reduced model; (c) speedup ratio

**Conclusions**

By integrating of the OSDCA for transient analysis with a POD-based approximated ESLs, a reduced order model (ROM) has been developed for dynamic topology optimization of the large-scale structures in the FE framework. In the suggested framework, the major contributions can be summarized as following:

(i) The OSDCA is proposed to establish the dynamic reduced order model in the dynamic finite element analysis stage;

(ii) The POD approximation strategy is employed to the displacement snapshots obtained by OSDCA method;

(iii) The OSDCA method integrates of the POD-based approximated ESLs, and it is applied to the topology optimization, thus an update objective function of dynamic topology optimization is defined. The efficiency of the proposed dynamic topology optimization strategy is reflected in that it reduces the complexity of the problem in both space and time domain. It not only reduces the DOF of the dynamic system, but also reduces the number of ESLs corresponding to time intervals, that is, it realizes the reduction in time domain. Consequently, the proposed reduction strategy has superior computational efficiency in dealing with large-scale problems.

From 2D and 3D numerical examples of different scales and different external loads, it can be concluded that the presented methodology has the ability to efficiently solve large-scale problems on the premise of ensuring the optimization results. Theoretically, the proposed OSDCA method might be used in different integral solution frameworks of motion equations. The current work also can be extended to the topology optimization of nonlinear large-scale dynamic systems.

**Declaration of competing interest**

The authors declare that they have no known competing financial interests or personal relationships that could have appeared to influence the work reported in this paper.

**Acknowledgements**

We acknowledge the support provided by the Project of the National Natural Science Foundation of China (11702090) and Peacock Program for Overseas High-Level Talents Introduction of Shenzhen City (KQTD20200820113110016).

**Replication of results**

A detailed description of the suggested dynamic topology optimization method has been provided in this paper. The authors are confident that the overall methodology can be reproduced. Due to the particularities of the project, the source codes and data cannot be shared. Readers interested in the data that support the findings of this study are encouraged to contact the corresponding author.